\documentclass[a4paper, 11pt]{article}
\usepackage[a4paper, margin=1in]{geometry}
\usepackage [utf8]{inputenc}
\usepackage{mathtools}
\usepackage{setspace}
\usepackage{graphicx}
\usepackage{float}
\usepackage{flexisym}
\usepackage{physics}
\usepackage{tabto}
\usepackage{amsfonts}
\usepackage{authblk}
\usepackage{enumitem}
\usepackage{amssymb}
\usepackage{amsthm}
\usepackage{amsmath}
\providecommand{\keywords}[1]{\textbf{\textit{Keywords---}} #1}

\theoremstyle{definition}

\title{{\Huge On Controllability of Driftless Control System on  Symmetric Spaces}}
\author[Archana Tiwari et al.]
{$^{1}$Archana Tiwari, $^{2}$R. N. Padhan  and $^{3}$K.C. Pati\footnote{Corresponding author:kcpati@nitrkl.ac.in}\\
{\small $^{1,3}$ Department of Mathematics, National Institute of Technology, Rourkela, Odisha-769008, India}\\
{\small $^{2}$Center for Applied Mathematics and Computing, ITER, Sikhsha 'O' Anusandhan, BBSR, Odisha-751030, India}\\

}
\date{}
\begin{document}
\maketitle
\begin{abstract}
\noindent {Symmetric spaces arise in wide variety of problems in Mathematics and Physics. They are mostly studied in Representation theory, Harmonic analysis and  Differential geometry. As many physical systems have symmetric spaces as their configuration spaces, the study of controllability on symmetric space is quite interesting.  In this paper, a driftless control system of  type $\dot{x}= \sum_{i=1}^m u_if_i(x)$ is considered on  a symmetric space. For this we have established global controllability condition which is illustrated by few examples of  exponential submanifolds of $SE(3)$ and random matrix ensembles.}
 
\end{abstract}


\section{Introduction}
Lie theory was first introduced by Brockett \cite{Bro1, Bro2} in motion control problems. He studied various aspects of controllability and observability of systems involving Lie groups.  The controllability properties of the systems on
 Lie groups were investigated in \cite{JS, JSK, SJ}.  Later, studies were done on these type of control systems from a geometric point of view \cite{Agr, Jean, Sach}.  Recent times have witnessed a lot of activities where by control theory (more particularly driftless control system) on Lie groups could play a lead role in explaining many physical, chemical and engineering problems in the area of quantum control \cite{Ram}, space altitude control \cite{Wal}, robotic controls \cite{Wu2}, chemical reaction control \cite{Dre, Far} etc. which could not be easily tackled by control theory defined on a conventional state space. With the advent of new technological challenges in control theory such as motion patterns of many kinesiological/robot mechanical systems or control on ensembles, it is being felt that a modified theory should be developed for tackling such type of situation, similar but not exact, to the case where Lie group replaced the ordinary state space in  control theory.
\par The motion pattern of human shoulder complex, human knee joint, mechanical device such as omni-wrist are associated with  spaces which can be approximated or shown to be some types of symmetric spaces. Similarly considering a dynamical system on random matrix ensembles such as Circular ensembles, Gaussian ensembles(Orthogonal (GOE), Sympletic (GSE), Unitary (GUE)) \cite{Kota} it becomes necessary to devise a theory of control system on symmetric space which is the main aim of our paper. We have confined ourselves to the study of controllability of non linear driftless control system on these type of spaces, as most of the physical system can be modelled on such background. However one can extend such types of studies to bilinear system or affine system etc. which we have not ventured into, in this paper.
\par Progresses have been made to some extent in this regard by various authors \cite{Sel, Kra} in which they have discussed optimal control problems in Lie groups where control  belong to symmetric space.  A class of optimal control problems defined on  certain kinds of symmetric spaces is discussed by Bloch \textit{et al.} \cite{Blo1, Blo2}. Similarly a lot of progress has been made to study control of complex networks \cite{Bag, Liu}  and control of ensembles  through structural  controllability \cite{Che1, Che2} and by other authors \cite{Agr1, Agr2, Lev, Li, Sam}, but to our knowledge nobody has discussed about these in the context of  random  matrix ensembles like GOE, GSE, GUE which actually correspond to Riemannian symmetric spaces classified by  Cartan \cite{Car}.  Taking these view points into consideration we were motivated to extend the study of control theory on Lie groups (Lie algebras) to symmetric spaces(Lie triple systems).
  \par Our task becomes simpler  as Lie triple systems and symmetric spaces are related by exponential mapping similar to Lie algebras and Lie groups. As a result many of the theories  on Lie groups are extended very easily to the case of symmetric spaces.  
\par In the next section we gave an introduction to the concepts and definitions of symmetric spaces and Lie triple system along with their relations. In the third section we define a nonlinear control system (driftless) on a symmetric space and obtain the condition under which a system is globally controllable. In the fourth section we apply these theories to a few simple but important examples. The last section contains few concluding remarks.

\section{Symmetric Space and Lie Triple System}
\subsection{Symmetric Space}
 \label{SS}
 A symmetric space\cite{Gil, Hel, Loos, Mun} is a manifold $S$ with a differential symmetric product $ \cdot $  that obeys the following conditions:
 \begin{enumerate}[label=(\roman*)]
 \item $a \cdot a=a$
 \item $a \cdot (a \cdot b)=b$
 \item $a \cdot (b \cdot c)= (a \cdot b)\cdot (a \cdot c)$,\\
 and moreover
 \item every $a$ has a neighbourhood $U$ such that  $ \forall b \in U$, $a \cdot b=b \Rightarrow b=a$.
 \end{enumerate}
 
 In other words, for any $x \in S$ there is some $s_x \in I(s)$  (the isometry group of $S$) with the properties
 \begin{equation*}
 s_x(x)=x, ~~~(ds_x)_x=-Id.
 \end{equation*}
This isometry $s_x$ is called symmetry at $x$. A symmetric space $S$ is precisely a homogeneous space with a symmetry $s_p$ at some point $p \in S $. A pointed symmetric space is a pair $(S,o)$ consisting of a symmetric space $S$ and a point $o$
called base point.
Some  examples of symmetric spaces are given below.
 \begin{enumerate}

\item The subsets of a  group $G$ which are closed under the composition $a \cdot b = ab^{-1}a$, where $ab$ is the usual multiplication in $G$. 

\item  The set of all symmetric positive definite matrices  forms a symmetric space with the product defined by $a \cdot  b = ab^{-1}a$.

\item Let $S=S^n$ be the unit sphere in $\mathbb{R}^{n+1}$ with the standard scalar product. The symmetry at any   $x \in S^n$  is the reflection at the line $\mathbb{R}x$ in $\mathbb{R}^{n+1}$, i.e., $s_x(y)=-y+ 2 \left\langle y,x  \right\rangle x.$ In this case, the symmetries generate the full isometry group which is the orthogonal group $O(n+1)$. 
\end{enumerate}
 \subsection{Lie Triple System}
Lie triple system \cite{Gil, Hel, Loos, Mun, Yad} is a vector space $L$ over a field $\mathbb{K}$ with a trilinear map $ \mu : L\times L\times L \longrightarrow L$ satisfying(if we write $\mu(A\otimes B \otimes C)=[A,B,C]$):
 \begin{enumerate}[label=(\roman*)]
\item $[A,A,B]=0$
\item $[A,B,C]+[B,C,A]+[C,A,B]=0$
\item $[A,B[D,E,F]]=[[A,B,D],E,F]+[D,[A,B,E],F]+[D,E,[A,B,F]]$, for $A,B,D,E,F \in L, [\cdot~,~\cdot,~\cdot]$ is called ternary operation of the Lie triple system.
\end{enumerate}

 In general, any subset of a Lie algebra $\mathfrak{g}$ that is closed under the operator 
\begin{center}
$T_A(\cdot)=ad^2_A(\cdot)=[A,[A,\cdot]]$
\end{center}
is a Lie triple system. 
\par  Some examples of Lie triple system is given below.
\begin{enumerate}

\item The set of all $n\times n$  matrices is a Lie triple system with Lie double bracket defined by $[A,B,C]=[[A,B],C]$, for any $n\times n$ matrices $A,B,C$  and  $[A,B]=AB-BA.$ In general, every Lie algebra $\mathfrak{g}$ with Lie product $[\cdot , \cdot]$ is a Lie triple system with ternary operation defined by $[ABC] = [[A,B],C]$, for any $A,B,C \in \mathfrak{g}$.

\item The set of all $n\times n$ symmetric matrices is a Lie triple system with Lie double bracket defined by $[A,B,C]=[[A,B],C]$, for $A,B,C$ any $n\times n$ symmetric matrices and  $[A,B]=AB-BA.$

\item Let $T$ be a Lie triple system and $S$ be any nonempty set. Then the set  $T^S$ of all functions $f$ from $S$ to $T$ is a Lie triple system with respect to a ternary operation given by $[f_1,f_2,f_3](x)=[f_1(x),f_2(x),f_3(x)].$

\end{enumerate}

 \subsection{Lie double bracket of vector fields on a symmetric space}

Let $X_1$ and $ X_2$ be smooth vector fields on  a $n$-manifold $M$, then $X_1$ and $ X_2$ are first order differential operators on smooth functions $M \longrightarrow \mathbb{R}$. Hence for a smooth function  $f: M \longrightarrow \mathbb{R}$,  $X_1f$ and $ X_2f$ are also smooth functions $M \longrightarrow \mathbb{R}$. 
\par Now we first define a differential operator $[X_1, X_2]$ termed as Lie bracket of $X_1$ and $ X_2$ as $[X_1, X_2]:=X_1X_2 - X_2X_1$, i.e., for  smooth functions  
 $f: M \longrightarrow \mathbb{R},$ $$[X_1, X_2](f):=X_1X_2(f) - X_2X_1(f).$$
 So by  definition $[X_1, X_2]$ is a second order differential operator. But it appears that  $[X_1, X_2]$ is indeed a first order differential operator and in fact, a vector field.
\par For any  smooth vector fields $X_1$ and $ X_2$, the Lie bracket $[X_1, X_2]$ is again a smooth vector field on $M$. Extending this, it can be shown that $[[X_1, X_2], X_3]$ is again a first order differential operator where $$[[X_1, X_2], X_3](f)=X_1X_2X_3(f) - X_3X_1X_2(f)- X_2X_1X_3(f)+X_3X_2X_1(f).$$

\par Before proceeding further, let us extend the  concepts of left  invariant vector fields on symmetric spaces. Let $\mathbf{T}$ be a Lie triple system $(LTS)$ with its ternary operation  $\mu(A\otimes B \otimes C)=[A,B,C]$ and $G$ be a Lie group. Then $G$ is said to act on $\mathbf{T}$ from left if there exist a function 
\begin{align*}
\varphi:& G \times \mathbf{T} \rightarrow \mathbf{T} \\
(g,A) & \mapsto  \varphi(g,A)=gA
\end{align*}
satisfying the following properties
\begin{enumerate}
\item $ex=x, ~ \forall x \in \mathbf{T}, $ where $e \in G$ is group identity,
\item $g_1(g_2x)=(g_1g_2)x, \forall  g_1, g_2 \in G ~\mbox{and}~ \forall x \in \mathbf{T}, $
\item for every $g\in G$, the left translation $\varphi_{g}=\varphi(g):\mathbf{T} \rightarrow \mathbf{T},~ A \mapsto gA$ is a linear map,
\item $\forall g \in G$ and $A,B,C, \in \mathbf{T};~ \mu(gA, gB, gC)=g\mu(A, B, C)=g[A,B,C].$
\end{enumerate}
We denote an action as above by $(G, \mathbf{T}).$ We call Lie triple system with an action of a group $G$ as $G-LTS.$
\par Further we see that the group $G$ acts on the space of points $G/K$ and maps it onto itself. In particular the coset representations $c \in G/K \subset G$ map the origin (base point) of $G/K$ on to the point $c \in G/K$. If $c$ is any point in $G/K$, it can be mapped into any other point $c'$ in $G/K$ by some group operation: $(c'c^{-1})c \rightarrow c', ~c'c^{-1}=g \in G.$  Moreover, the only group operation which leaves every point $c$ fixed is the identity, $gc=c, ~\forall c \in G/K \Rightarrow g=Id.$
\par The vector fields  of the form $U(X)=XA,~ V(X)=XB,~ W(X)=XC,~~X \in S,~~A,B,C \in \mathbf{T}$ (tangent space at the base point)  are called left invariant vector fields if they satisfy the following property
\begin{equation*}
[[U,V],W](X)=[XA,XB,XC]=X[A,B,C].
\end{equation*}

\par Suppose now we have three left invariant  smooth vector fields  $U,V,W$  on a symmetric space $S$ and  $V_F(S)$ denotes the space of all smooth vector fields on $S$. The Lie double bracket of the fields $U,V,W$ is the vector field, $[[U,V],W]] \in V_F(S).$ It can be shown that
\begin{center}
$[[U,V],W](X)= \frac{d}{dt}\bigg|_{t=0} \gamma (\sqrt[3]{t}),~~~~~X \in S,$
\end{center}
where the curve $\gamma$ is defined as 
\begin{equation}
\gamma(t)=e^{-tW} \circ  e^{-tU} \circ e^{-tV} \circ e^{tU} \circ e^{tV} \circ e^{tW} \circ e^{-tV} \circ e^{-tU} \circ e^{tV} \circ e^{tU}(X).
\label{db}
\end{equation}


Here $e^{tU}$ denotes the flow of the vector field $U$
\begin{center}
$ \frac{d}{dt}e^{tU}(X)=U(e^{tU}(X)),~~~~e^{tU} \bigg|_{t=0}(X)=X. $
\end{center}

Now, we have  the  Lie double bracket of the left invariant vector fields 
\begin{align*}
[[U,V],W](X)&=[[XA,XB],XC]\\
&=X[[A,B],C]\\
&=X((AB-BA)C+C(BA-AB))\\
&=X(ABC-BAC+CBA-CAB).
\end{align*}
\par As flows of the left invariant vector fields are given by the matrix exponential:
\begin{center}
$e^{tU}(X)=X e^{(tA)},~~~~e^{tV}(X)=X e^{(tB)},~~~~e^{tW}(X)=X e^{(tC)}$.
\end{center}
Computing the lower order terms of the curve $\gamma$ from  (\ref{db}); we obtain
\begin{align*}
\gamma(t)=& X e^{(tA)}e^{(tB)}e^{(-tA)}e^{(-tB)}e^{(tC)}e^{(tB)}e^{(tA)}e^{(-tB)}e^{(-tA)}e^{(-tC)}\\
=& X \Bigg(Id+tA+ \frac{t^2}{2!}A^2+ \frac{t^3}{3!}A^3+\cdots \Bigg) \Bigg(Id+tB+ \frac{t^2}{2!}B^2+ \frac{t^3}{3!}B^3+ \cdots \Bigg)
\Bigg(Id-tA+ \frac{t^2}{2!}A^2- \frac{t^3}{3!}A^3+\cdots \Bigg)\\& \Bigg(Id-tB+ \frac{t^2}{2!}B^2- \frac{t^3}{3!}B^3+\cdots\Bigg) \Bigg(Id+tC+\frac{t^2}{2!}C^2+ \frac{t^3}{3!}C^3+\cdots \Bigg)\Bigg(Id+tB+ \frac{t^2}{2!}B^2+ \frac{t^3}{3!}B^3+\cdots \Bigg)\\ & \Bigg(Id+tA+ \frac{t^2}{2!}A^2+ \frac{t^3}{3!}A^3+\cdots \Bigg)\Bigg(Id-tB+ \frac{t^2}{2!}B^2- \frac{t^3}{3!}B^3+\cdots \Bigg) \Bigg(Id-tA+ \frac{t^2}{2!}A^2- \frac{t^3}{3!}A^3+\cdots \Bigg) \\ & \Bigg(Id-tC+ \frac{t^2}{2!}C^2- \frac{t^3}{3!}C^3+\cdots \Bigg)\\
=&X \Bigg(Id+t^2(AB-BA)+ \frac{t^3}{2!}(A^2B+BA^2+2BAB-2ABA-AB^2-B^2A)\Bigg)\\ &\Bigg(Id+tC+ \frac{t^2}{2!}C^2+ \frac{t^3}{3!}C^3+ \cdots \Bigg)\\ & ~\Bigg (Id+t^2(BA-AB)+ \frac{t^3}{2!}(AB^2+B^2A+2ABA-2BAB-A^2B-BA^2)\Bigg)\\ &\Bigg(Id-tC+ \frac{t^2}{2!}C^2- \frac{t^3}{3!}C^3+\ldots\Bigg) \\
=&X \Bigg(Id +tC+ \frac{t^2}{2!}(C^2+2AB-2BA)+ \frac{t^3}{3!}(C^3+6ABC-6BAC+3A^2B+3BA^2\\ &+6BAB-6ABA-3AB^2-3B^2A)\Bigg)  \Bigg(Id-tC+ \frac{t^2}{2!}(C^2+2BA-2AB)\\ &- \frac{t^3}{3!}(C^3+6BAC-6ABC-3A^2B-3B^2A+6BAB-6ABA+3A^2B+3BA^2)\Bigg)\\
=&X \Bigg(Id+t^3(AB-BA)C-C(AB-BA)\Bigg)\\
=&X \Bigg(Id+t^3[[A,B],C]+\cdots \Bigg).\\
\end{align*}
Hence $$\gamma (\sqrt[3]{t})=X(Id+t[[A,B]C]+ \cdots ),$$

 is a smooth curve at $t=0$, and 
\begin{center}
$\frac{d}{dt}\bigg|_{t=0} \gamma (\sqrt[3]{t})=[[A,B],C].$
\end{center}
Thus, we see that these types of vector fields form a tangent space to the symmetric space at the base point and is called Lie triple system.

\par Special types of symmetric spaces can also be constructed as follows:\\
A symmetric space is associated to an involutive automorphism of a given Lie algebra $\mathfrak{g}$. To be more specific, if $\sigma$ is an automorphism it preserves multiplication: $ [\sigma(x), \sigma(y)]=\sigma([x,y]).$ Suppose that the linear automorphism $\sigma:\mathfrak{g} \rightarrow \mathfrak{g}$ is such that $\sigma^2=1,$ but $\sigma$  is not the identity, that means $\sigma$ has eigen values $\pm 1$ and it splits the algebra $\mathfrak{g}$ into orthogonal eigen spaces corresponding to these eigen values. Such a mapping is called an involutive automorphism.
\par Suppose $\mathfrak{g}$ is a compact simple Lie algebra, $\sigma$ is an involutive automorphism of $\mathfrak{g}$ and $\mathfrak{g}=k\bigoplus p$ where 
\begin{align*}
\sigma(x)&=x, ~~~\mbox{for}~x \in k,\\
\sigma(x)&=-x, ~~~\mbox{for}~x \in p.
\end{align*}

It is easy to check that $k$ is a subalgebra but $p$ is not. In fact, the commutation relation 
\begin{equation}
[k,k] \subset k, ~~~[k,p] \subset k,~~~[p,p]\subset k
\label{cr}
\end{equation}
holds. A subalgebra `$k$' satisfying (\ref{cr}) is called a symmetric subalgebra. If we now multiply the element of $p$ by $i$ (weyl unitary trick), we construct a new non compact algebra $\mathfrak{g}^{*}=k+p^{*}$. This is called Cartan decomposition and `$k$' is the maximal compact subalgebra of $\mathfrak{g}^{*}$. The coset spaces $P\simeq e^p \simeq G/K$, where $G\simeq e^{\mathfrak{g}}$, is a symmetric space. Similarly $G^{*}/K$ is also a symmetric space. The corresponding Lie triple systems are denoted by $\mathfrak{g}/k$ and $\mathfrak{g}^{*}/k$ respectively. For example
\begin{align*}
G/K &=SU(n, \mathbb{C})/SO(n, \mathbb{R});\\
G^{*}/K &=SL(n, \mathbb{R})/SO(n, \mathbb{R});
\end{align*}
are symmetric spaces of compact and non compact type respectively. These types of symmetric spaces has been classified by Cartan  \cite{Car} and they corresponds to various random matrix ensembles. 
\par The Lie triple systems associated with the above two types of symmetric spaces are denoted by 
\begin{align*}
\mathfrak{g}/k &=su(n, \mathbb{C})/so(n, \mathbb{R});\\
\mathfrak{g^{*}}/k &=sl(n, \mathbb{R})/so(n, \mathbb{R});
\end{align*}
respectively.
 \section{Driftless Control System on a Symmetric space and its Controllability Condition} 
 
 To begin with, we introduce some basic concepts and definitions related with control system following standard literature on control theory.
 
 \par Let $S$ be a smooth  $n$-dimensional manifold with tangent space at $x$ denoted by $T_xS$. A general control system takes the form

 $$\dot{x}=f(x,u)$$ 
 where $x \in S$ denotes the state, $u \in \mathbb{R}$ is the control,  and $f(\cdot, u), $ is a vector field on $S$ for all  $ u $.
 \par Suppose that $f$ is locally Lipschitz relative to the second variable then  $ \forall x\in S ~\mbox{and}~ \forall u \in L^2([0, t_f], \mathbb{R}^m)$, the cauchy problem
 \begin{equation*}
 \dot{x}(t)=f(x(t), u(t)), t \in [0, t_f],~ x(0)=x
  \end{equation*}
has a unique solution 
    \begin{equation*}
 x(\cdot)=x(\cdot, x, u):[0, t_f']\longrightarrow S ~\mbox{with}~ t_f' \leq t_f.
  \end{equation*}

 Mainly there are three  types of controllability: global controllability, local controllability at an equilibrium point and local controllability along a reference trajectory. In this paper we extend the study of global controllability to symmetric spaces.

 \par A  control system of the form $\dot{x}=f(x,u)$ defined as above is said to be globally controllable on $S$ if for any $x_1, x_2 \in S $ and $t_f > 0$ there exist a control $u \in L^2([0, t_f], \mathbb{R}^m)$ such that the solution of the cauchy problem starting at $x_1$ satisfy $x(t_f)=x_2$.

Now, given a family $\mathcal{F}$ of smooth vector fields on  symmetric space $S$, denote $LTS \left\lbrace \mathcal{F} \right\rbrace $ as the Lie triple system generated by $\mathcal{F}$. It is the smallest vector subspace $\mathcal{V}$ of smooth vector fields containing $\mathcal{F}$, which  satisfies  
  
  \begin{equation*}
  [[A,B],C] \in \mathcal{V},~ \forall A,B \in \mathcal{F}, ~\forall C \in \mathcal{V}.
  \end{equation*}
   
  The sufficient condition of global controllability for a driftless control system on Lie group is given and proved independently by P. Rashesvsky \cite{Ras} and W. L. Chow \cite{Cho}. More recently simpler proofs given by F.  Jean \cite{Jean} and L. Rifford \cite{Rif} separately. In this paper we have extended Rashesvsky and Chow's theorem to symmetric spaces following Jean closely, which we call extended Rashesvsky and Chow's theorem for symmetric spaces, which is explained below:
  
\textbf{Extended Rashesvsky-Chow's Theorem.} Let $S$ be a smooth connected $n$-manifold of symmetric space and $X_1, X_2, \ldots, X_m $ be $m$ smooth vector fields on $S$. Assume that $$LTS \left\lbrace X_1, X_2, \ldots, X_m \right\rbrace (x)=T_xS, ~ ~ \forall x \in S .$$
 Then the control system 
  \begin{equation*}
   \dot{x}= \sum_{i=1}^m u_iX_i(x), 
  \end{equation*}
  
is globally controllable on $S$.

  The controllability of a driftless control system on a symmetric space is characterized by the properties of the Lie triple system generated by the vector fields $\left\lbrace X_1, X_2, \ldots, X_m\right\rbrace . $
  \par Let $V_F(S)$ denote the set of all smooth vector fields on $S$ and  $\Theta^1$  be the linear subspace of $V_F(S)$ generated by vectors fields $\left\lbrace X_1,\cdots, X_m\right\rbrace $,
\begin{center}
i.e., $\Theta^1=span \left\lbrace X_1,\cdots, X_m\right\rbrace. $
\end{center}
For $p \geq 1$, define 
\begin{center}
$\Theta^{p+1}=\Theta^p + [ [\Theta^1, \Theta^p],\Theta^1]$
\end{center}
where  $[ [\Theta^1, \Theta^p], \Theta^1]=span \left\lbrace[[A,B],C]]:A, C \in \Theta^1, B \in \Theta^p\right\rbrace $.\\
The Lie triple system generated by $\left\lbrace X_1,\cdots, X_m \right\rbrace $ is defined to be
\begin{center}
$LTS \left\lbrace X_1,\cdots, X_m \right\rbrace = \bigcup_{p \geq 1} \Theta^p.$
\end{center}

 From the  Jacobi identity property   of Lie triple system, $LTS\left\lbrace X_1,\cdots, X_m \right\rbrace $ is the smallest linear subspace of $V_F(S)$ which  contains $X_1,\cdots, X_m$ and is invariant under Lie double brackets.
\par  Let  $\mathcal{L}  $ denote the free Lie algebra generated by the elements $\left\lbrace  1, \ldots, m \right\rbrace. $ Here,   $\mathcal{L} $ is the $\mathbb{R}$ vector space generated by $\left\lbrace  1, \ldots, m \right\rbrace $ and their bracket $[\cdot , \cdot]$ combined with relation of skew symmetry and the Jacobi identity.
\par The length of an element $I$ of $\mathcal{L}$ is denoted by $|I|$ and is defined inductively by $|I|=1$ for $I=1,\ldots, m$, $|I|=|I_1|+|I_2|$ for $I= [I_1, I_2]$ and $|I|=|I_1|+|I_2|+|I_3|$ for $I= [I_1, I_2], I_3].$ With every element $I \in \mathcal{L} $  associate the vector field $X_I \in LTS \left\lbrace  X_1,\cdots, X_m \right\rbrace $ obtained by plugging in $X_i,i=1,2,\ldots,m$ for the corresponding letter $i$ in $I$. For example $X_{[1,2,3]}=[[X_1, X_2],X_3]. $
\par Due to  Jacobi identity property  for Lie triple system
\begin{equation*}
\Theta^p=span \left\lbrace X_I:|I|\leqslant p \right\rbrace,
\end{equation*}.

For $x \in S$,  we set $LTS \left\lbrace X_1,\cdots, X_s \right\rbrace  (x)=\left\lbrace X(x): X \in LTS(X_1,\cdots, X_m) \right\rbrace $ and for $p \geqslant 1, \Theta^p(x)=\left\lbrace X(x): X \in \Theta^p \right\rbrace$. By definition, these sets are linear subspaces of $T_xS$.

  \par We say that the control system   
  \begin{equation}
   \dot{x}= \sum_{i=1}^m u_iX_i(x),  
   \label{eqn:driftless}
  \end{equation}
  and the vector fields $(X_1,\cdots, X_m)$ satisfies Chow's condition if $LTS(X_1,\cdots, X_m)(x)=T_xS, ~\forall x \in S.$

Equivalently, for any $x \in S ~\exists $ an integer $r=r(x)$ such that $dim \Theta^r(x)=n.$ 

\par If \eqref{eqn:driftless} satisfies Chow's condition, then  $\forall x \in S$, the reachable set  $R_x$ is a neighbourhood of $x$.

\par Now we confine ourselves to a  small neighbourhood of $x$, i.e., $U \subset S$ which can be identified with a neighbourhood of $0$ in $\mathbb{R}^n.$\\

Let $\Psi^i_t=exp(tX_i), i=1,2,\cdots, m$ be the flow of the vector field $X_i.$ Each curve $t \longmapsto  \Psi^i_t(x_1)$ is a trajectory of \eqref{eqn:driftless} and we have
\begin{center}
$\Psi^i_t=id+tX_i+O(t)$.
\end{center}

For every element $I \in \mathcal{L} $, we define the local diffeomorphism $\Psi^I_t$ on $U$ by induction on the length $|I|$ of $I:$ if $I=[I_1, I_2], I_3]$  then
\begin{center}
$\Psi^I_t= [[\Psi^{I_1}_t,\Psi^{I_2}_t],\Psi^{I_3}_t]= \Psi^{I_3}_{-t} \circ \Psi^{I_1}_{-t}\circ \Psi^{I_2}_{-t}\circ \Psi^{I_1}_t \circ \Psi^{I_2}_{t}\circ \Psi^{I_3}_{t} \circ \Psi^{I_2}_{-t}\circ \Psi^{I_1}_{-t}\circ \Psi^{I_2}_t \circ \Psi^{I_1}_t. $
\end{center}
By construction, $\Psi^I_t$ may be expanded as a composition of flows of the vector fields $X_i, i=1,2,\cdots, m.$  As a result $\Psi^I_t(x_1)$ is the end point of a trajectory of \eqref{eqn:driftless} induced from $x_1$. Further, on a neighbourhood of $x$ there holds
\begin{equation}
\Psi^I_t=id+t^{|I|}X_I+O(t^{|I|}).
\end{equation}

A  diffeomorphism now can be defined whose derivative with respect to the time is exactly $X_I$ and is given by 
\begin{equation}
\psi^I_t=id+tX_I+O(t)
\label{cont4}
\end{equation}
and $\psi^I_t(x_1)$ is the endpoint of a trajectory of \eqref{eqn:driftless} starting from $x_1$.
\par Now we choose double commutators $X_{I_1},X_{I_2},\cdots, X_{I_n}$ whose values at $x$ span $T_{x}S$. This is feasible through Chow's condition. Let $\phi$ be a map defined on a small neighbourhood $U_0$ of $0$ in $\mathbb{R}^n$ by
\begin{center}
$\phi(t_1,\cdots, t_n)=\psi^{I_n}_{t_n} \circ \cdots \circ \psi^{I_1}_{t_1}(x_1) \in S.$
\end{center}
We infer from (\ref{cont4}) that this map is $C^1$ near 0 and has an invertible derivative at 0, which implies that it is a local $C^1$-diffeomorphism. Hence $\phi(U_0)$ contains a neighbourhood of $x$.

\par Now,  $ \forall t \in U_0, ~\phi(t)$ is the endpoint of a concatenation of trajectories of \eqref{eqn:driftless}, the first one starting from $x$. It is then the endpoint of a trajectory starting from $x$. Therefore $\phi(U_0) \subset R_{x}$, which implies that $R_{x}$ is a neighbourhood of $x$.

Let $x_1 \in S$ and if $x_2 \in R_{x_1}$, then $x_1 \in R_{x_2}.$ As a result, $R_{x_1}=R_{x_2}$ for any ${x_2} \in S$ and by the above result, $R_{x_1}$ is an open set. Therefore the union of the sets $R_{x_1}$ which  are pairwise disjointed covers the  manifold $S$. Since $S$ is a connected manifold, there can only be one open set which can cover $S$, i.e.,  $R_{x_1}$. Hence any two points in $S$ can be connected by trajectories of \eqref{eqn:driftless}.

\par Thus we see that as \eqref{eqn:driftless} satisfies the Chow's condition, $\forall x \in S$, the reachable set $R_x$ is a neighbourhood of $x$. So now if $S$ is connected and if  \eqref{eqn:driftless} satisfies the Rashesvsky-Chow's condition then any two points of $S$ can be connected by a trajectory of \eqref{eqn:driftless} and we obtain the extended Chow's theorem.

\section{Applications}
Immediate examples of symmetric spaces which are applied directly to physical applications are symmetric submanifolds of special Euclidean group $SE(3)$ which are related to various  kinesiological and mechanical systems and accordingly, have numerous potential applications in robot kinematics \cite{Wu1, Wu2}. Similarly,  as discussed earlier the random matrix ensembles like GOE, GUE and circular ensembles have a lot application in quantum transport problems, etc. So in this section  we have studied  the controllability aspect of some symmetric spaces.  But to start with, we have simplest example of symmetric space $SO(3)/SO(2)\simeq S^2$.  Each element of $SO(3)/SO(2)$ can be fully characterized by three real parameters such that their moduli sum to 1 and then there is a one-to-one correspondence between each element of $SO(3)/SO(2)$ and a set of Cartesian coordinates for $S^2$.  

\subsection{Controllability on $SO(3)/SO(2)$}
The infinitesimal generators of the Lie algebra $so(3)$ of Lie group $SO(3)$, corresponds to the derivative of rotation around the each of the standard axes, evaluated at the identity are:
\[
X_1=
\begin{bmatrix}
     0 & 0  & 0\\
    0 & 0  & -1\\
    0 & 1  & 0
\end{bmatrix},
X_2=
\begin{bmatrix}
    0 & 0 & 1\\
    0 & 0 & 0 \\
    -1 & 0 & 0 
\end{bmatrix},
X_3=
\begin{bmatrix}
    0 & -1  & 0\\
    1 & 0  & 0\\
    0 & 0  & 0 
\end{bmatrix}.
\]
The infinitesimal generators of the Lie algebra $so(2)$ of Lie group $SO(2)$ correspond to the derivative of rotation evaluated at the identity is:
\[
X=
\begin{bmatrix}
     
     0  & -1\\
    1 	& 0 
\end{bmatrix},
\]

This result for infinitesimal generator for rotation about the $z$- axis i.e. $X_3$ is essentially identical to those of $SO(2)$. So the generators of $so(3)/so(2)$ are $ \lbrace X_1,X_2 \rbrace$ and 
\[
\left[  X_1,X_2\right]  = 
\begin{bmatrix}
    0 & -1  & 0\\
    1 & 0  & 0\\
    0 & 0  & 0 
\end{bmatrix}=X_3.
\]
Also $\left[   \left[  X_1,X_2\right], X_1 \right] = \left[  X_3,X_1\right]= X_2$ and $\left[   \left[  X_1,X_2\right], X_2 \right] = \left[  X_3,X_2\right]= -X_1$. We can verify that $ \left\lbrace  X_1, X_2 \right\rbrace $ form a Lie triple system. 

Now we consider a  driftless control system of the following type in classical notation \cite{Sach} 
\begin{equation}
 \label{SO}
  \dot{x}= x\sum_{i=1}^2 X_iu_i, ~x \in S
  \end{equation}
where $x \in SO(3)/SO(2), ~u_i $ are controls on it and $X_i$ are smooth vector fields on the given symmetric space.

Let $\Theta^1= span \lbrace X_1,X_2 \rbrace$,
for $p>1, \Theta^{p+1}= \Theta^p + \left[ \Theta^1, \left[ \Theta^1, \Theta^p\right] \right] $,
so $ \Theta^2= \Theta^1 + \left[ \Theta^1, \left[ \Theta^1, \Theta^1\right] \right]$\\
Let $A, B, C \in \Theta^1$ be arbitrary, then we can write\\ $A= \alpha_1X_1+\beta_1X_2,~~~ B= \alpha_2X_1+\beta_2X_2,~~~ C= \alpha_3X_1+\beta_3X_2,~~~~~\alpha_i, \beta_i, i\in \left\lbrace 1,2,3\right\rbrace $ are scalars. \\
Now, $\left[ \Theta^1, \left[ \Theta^1, \Theta^1\right] \right]= \left[ \alpha_1X_1+\beta_1X_2, \left[ \alpha_2X_1+\beta_2X_2, \alpha_3X_1+\beta_3X_2\right] \right]$\\
$=\left[ \alpha_1X_1+\beta_1X_2, \left(  \left[ \alpha_2X_1, \alpha_3X_1 \right] + \left[ \alpha_2X_1, \beta_3X_2 \right] + \left[  \beta_2X_2, \alpha_3X_1\right] + \left[  \beta_2X_2, \beta_3X_2 \right] \right)  \right]$ \\
$=\left[ \alpha_1X_1 +\beta_1X_2, \alpha_2\beta_3\right[ X_1,X_2 \left]  + \beta_2\alpha_3 \right[ X_2,X_1 \left]\right] $\\
$=\left[ \alpha_1X_1+ \beta_1X_2,\left(  \alpha_2\beta_3-\beta_2\alpha_3 \right) X_3\right] $\\
$=\left[ \alpha_1X_1+ \beta_1X_2, \gamma X_3\right],~ \mbox{where} ~\gamma= \alpha_2\beta_3-\beta_2\alpha_3 $\\
$=\alpha_1 \gamma \left[ X_1,X_3\right]  +  \beta_1 \gamma \left[  X_2,X_3 \right] $\\
$=\alpha_1 \gamma X_1 + \beta_1 \gamma  X_2 \in \Theta^1$\\
So, $\Theta^2= \Theta^1 + \left[ \Theta^1, \left[ \Theta^1, \Theta^1\right] \right]$\\
$=\Theta^1$\\
Hence, number of generators of  $\Theta^2=2$. This implies dim $\Theta^2=2$.\\
And hence $so(3)/so(2)$ satisfies the Chow's condition as  dim $\Theta^2=2$ i.e., dim $\Theta^2$ is equal to the dimension of the symmetric space $SO(3)/SO(2)$. So we conclude that the control system defined in (\ref{SO}) is controllable by  Chow-Rashevsky's theorem. This is a trivial case, however we have done it in detail and the result coincides with  the paper by Brockett for controllability on $S^n$ sphere \cite{Bro2, Lob, Sach}. 

\subsection{Controllability on symmetric submanifolds of  SE(3)}
The special Euclidean group $SE(3)$ admits an inversion symmetry through any of its elements hence it is a symmetric space.  The Lie algebra of $SE(3)$ is $se(3)$  whose basis elements are these $4\times 4$ matrices, each of which corresponds to either infinitesimal rotations or infinitesimal translations along each axis:

\[
e_1=
\begin{bmatrix}
    0 & 0 & 0 & 1\\
    0 & 0 & 0 & 0\\
    0 & 0 & 0 & 0\\
    0 & 0 & 0 & 0
\end{bmatrix},
e_2=
\begin{bmatrix}
    0 & 0 & 0 & 0\\
    0 & 0 & 0 & 1\\
    0 & 0 & 0 & 0\\
    0 & 0 & 0 & 0
\end{bmatrix},
e_3=
\begin{bmatrix}
    0 & 0 & 0 & 0\\
    0 & 0 & 0 & 0\\
    0 & 0 & 0 & 1\\
    0 & 0 & 0 & 0
\end{bmatrix},
\]
\[
e_4=
\begin{bmatrix}
    0 & 0 & 0 & 0\\
    0 & 0 & -1 & 0\\
    0 & 1 & 0 & 0\\
    0 & 0 & 0 & 0
\end{bmatrix},
e_5=
\begin{bmatrix}
    0 & 0 & 1 & 0\\
    0 & 0 & 0 & 0\\
    -1 & 0 & 0 & 0\\
    0 & 0 & 0 & 0
\end{bmatrix},
e_6=
\begin{bmatrix}
    0 & -1 & 0 & 1\\
    1 & 0 & 0 & 0\\
    0 & 0 & 0 & 0\\
    0 & 0 & 0 & 0
\end{bmatrix}.
\]

There is one to one correspondence between Lie triple systems of $SE(3)$ and the symmetric submanifolds of $SE(3)$. Therefore, the classification of symmetric submanifolds of $SE(3)$ upto conjugation is equivalent to that of the conjugacy class of Lie triple system of $SE(3)$. There are seven conjugacy classes of symmetric submanifolds of $SE(3)$, all of which can be locally represented  by $ exp~M$ with $M$ being a Lie triple (sub)system of $se(3)$.  As discussed earlier, the Lie triple system are  vector subspace of $se(3)$ closed under double Lie brackets. The Lie triple systems of $se(3)$ \cite{Wu2} are $ \lbrace e_3, e_4 \rbrace, ~ \lbrace e_3, e_4+ \mathsf{p}e_1 \rbrace, ~ \lbrace e_4, e_5 \rbrace, ~ \lbrace e_1, e_3, e_4 \rbrace, ~ \lbrace e_3, e_4, e_5 \rbrace,~  \lbrace e_1, e_2, e_4, e_5 \rbrace $ and $\lbrace e_1, e_2, e_3, e_4, e_5 \rbrace $ here $\mathsf{p}$ is the pitch value  of the generic screw which takes an arbitrary finite value.\\
 From control theory point of view, a system is controllable  on a space with all available controls in hand. But generally we are interested in controlling a system with fewer number of controls. In this case we see that the following two systems are controllable on the given symmetric spaces.
\begin{enumerate}

\item A  left invariant driftless control system defined on the symmetric submanifold of $SE(3)$, represented by  $M_1=exp \left\lbrace e_1, e_2, e_4, e_5\right\rbrace$ is given by 
\begin{equation}
 \label{se1}
  \dot{x}= x\left(  e_1u_1+ e_2u_2+e_4u_4+ e_5u_5\right), ~x \in M_1, ~u_i, i=1,2, 4,5 \in \mathbb{R}.
  \end{equation}

Since, $ \left[ \left[ e_2,e_4 \right] , e_5\right] =e_1$, so the three generators $\left\lbrace e_2, e_4, e_5  \right\rbrace $ upon double bracket commutation generates the full set of basis elements of $LTS\left\lbrace e_1, e_2, e_4, e_5\right\rbrace$. So the system (\ref{se1}) can be  controlled with fewer number of controls by  Rashvesky-Chow's theorem and  is represented as
\begin{equation}
 \label{se1a}
  \dot{x}= x\left(  e_2u_2+e_4u_4+ e_5u_5+\right), ~x \in M_1, ~u_i, i=2,4,5 \in \mathbb{R}.
  \end{equation}

\item Similarly, a left invariant driftless control system defined on the symmetric submanifold of $SE(3)$, represented by  $M_2=exp \left\lbrace e_1, e_2,e_3,  e_4, e_5\right\rbrace$ is given by 
\begin{equation}
 \label{se2}
  \dot{x}= x\left(  e_1u_1+ e_2u_2+e_3u_3+e_4u_4+ e_5u_5\right), ~x \in M_2, ~u_i, i=1,\ldots,5 \in \mathbb{R}.
  \end{equation}

Since, $ \left[ \left[ e_1,e_5 \right] , e_4\right] =e_2$, so, the four generators $\left\lbrace e_1, e_3,e_4, e_5  \right\rbrace $ upon double bracket commutation generates the full set of basis elements of $LTS\left\lbrace e_1, e_2, e_3, e_4, e_5\right\rbrace$. So the system (\ref{se2}) can be rewritten as 
\begin{equation}
 \label{se2a}
  \dot{x}= x\left(  e_1u_1+e_3u_3+e_4u_4+ e_5u_5\right), ~x \in M_2, ~u_i, i=1,3,4,5 \in \mathbb{R}.
  \end{equation}

and it is controllable by  Rashvesky-Chow's theorem.

\end{enumerate}

\subsection{Controllability on symmetric spaces associated with random matrices}
Random matrix theory deals with the study of matrix ensembles, i.e., matrices with a probability measure. The classical non-compact matrix Gaussian ensembles are  GOE, GUE and GSE. These  ensembles arises when a Gaussian-like probability measure is introduced in each set of a family of non-compact matrices.  Each  set of matrices, along with the probability measure on it, is invariant under the action of the orthogonal, unitary or symplectic group, respectively. The ensembles above have compact counterparts  the Circular Orthogonal Ensemble (COE), the Circular Unitary Ensemble (CUE) and the Circular Symplectic Ensemble (CSE),  leaving the probability measure invariant. The members of all the six ensembles discussed  are submanifolds of Euclidean space. In fact, they are (essentially) Riemannian globally symmetric spaces. The integration manifolds of the random matrices, the distribution of eigenvalues and the Dyson  characterizing the ensembles \cite{Haa} are strictly in correspondence with symmetric spaces.  Many important results can be acquired from this identification. A new classification of random matrix ensembles arises from the Cartan classification of triplets of symmetric spaces with positive, zero and negative curvature. Now we discuss controllability on two examples of GOE and COE respectively.

\subsubsection{Gaussian Orthogonal Ensemble (GOE) }
The GOE is the ensemble formed, for each $n$, by the set of $n\times n$ real symmetric matrices, whose entries are independent normal random variables
\begin{equation}
S_{GOE}(n)=\left\lbrace A=(a_{ij})_{n \times n}=(a_{ji})\right\rbrace
\label{goe}
\end{equation}
 with the probability measure
 \begin{equation}
 d\mu(A) \propto exp(-tr \frac{A^2}{2})dA
 \end{equation}
 where $ dA=\prod_{i\leqslant j}da_{ij}$ is the (additive) Haar measure on $S_{GOE}$. The nonsingular matrices in (\ref{goe}) form an open set of full measure, one of whose components $S$ consists of all the real symmetric positive definite $n\times n$ matrices. It is  known the symmetric positive definite matrices \cite{Moa} are symmetric spaces, with symmetric product $a\cdot b=ab^{-1}a$ and its Lie triple system are the set of symmetric matrices.  The Lie triple system of the $3\times 3$ symmetric positive random matrix, are the set of  $3\times 3$ symmetric matrices, whose basis elements are 
  
\[
a_1=
\begin{bmatrix}
    1 & 0 & 0 \\
    0 & 0 & 0 \\
    0 & 0 & 0 
\end{bmatrix},
a_2=
\begin{bmatrix}
    0 & 0 & 0 \\
    0 & 1 & 0 \\
    0 & 0 & 0 
\end{bmatrix},
a_3=
\begin{bmatrix}
    0 & 0 & 0 \\
    0 & 0 & 0 \\
    0 & 0 & 1 
\end{bmatrix},
\]
\[
a_4=
\begin{bmatrix}
    0 & 1 & 0 \\
    1 & 0  & 0\\
    0 & 0 & 0
\end{bmatrix},
a_5=
\begin{bmatrix}
    0 & 0 & 1 \\
    0 & 0 & 0 \\
    1 & 0 & 0 
\end{bmatrix},
a_6=
\begin{bmatrix}
    0 & 0 & 0 \\
    1 & 0 & 1 \\
    0 & 1 & 0 

\end{bmatrix}.
\]
So a driftless control system on the symmetric space $S$ of symmetric  positive definite matrices with six number of controls  is given by 
\begin{equation}
 \label{random1}
  \dot{x}= x\sum_{i=1}^6 a_iu_i, ~x \in S, ~u_i \in \mathbb{R}. 
  \end{equation}

Since, $[[a_1,a_4], a_4]=-2a_2, [[a_2,a_4], a_4]=a_6 ~\mbox{and}~ [[a_2,a_6], a_4]=a_5$. The three basis element $\left\lbrace a_1, a_3, a_4 \right\rbrace $ upon double bracket commutation generates the full set of basis elements of the $3 \times 3$ symmetric matrices. So now the control system with fewer number of controls is given by 
\begin{equation}
 \label{random1a}
  \dot{x}= x\left(  a_1u_1+a_3u_3+a_4u_4\right) , ~x \in S, ~u_i \in \mathbb{R}. 
  \end{equation}
  
It is controllable on given GOE by Rashvesky-Chow's theorem. 
 The maximum number of controls is $6$. We can construct a number of control systems with different combinations  with fewer number of controls, but not all such systems are controllable. In particular if we take the control system (\ref{random1a}) as above then we see that only  with these three number of controls, we can control the whole system defined in (\ref{random1}).
\subsubsection{ Circular Orthogonal Ensemble (COE)}
 In $SU(n)$, the subgroup $SO(n)$ is the fixed-point set of the involution 
\begin{center}
 $a\mapsto a^{\sigma}:=(a^{-1})^t$.
\end{center}
The space $SU(n)/SO(n)$ can be realized as the set of matrices
\begin{equation}
S(n)=\left\lbrace A \in SU(n)|~ A ~\mbox{is symmetric, i.e.,}~ A=A^T\right\rbrace.
\end{equation}
The ensemble $S(n)$ endowed with its Haar probability measure is Dyson's Circular Orthogonal Ensemble (COE) \cite{Case} i.e., the integration manifold of the circular orthogonal ensembles  is $SU(n)/SO(n)$. If we have consider the case of $n=3$, 
the group $SU(3)$  is characterized by $3\times 3$ unitary  matrices with unit determinant. The generators are Hermitian
matrices $3 \times 3$ with zero trace. Such general Hermitian matrix can be parametrized with eight real numbers $a, \ldots,  h$
\[
\begin{bmatrix}
     a & c-id  & e-if\\
    c+id & b  & g-ih\\
    e+if & g+ih  & -a-b
\end{bmatrix}.
\]
The Lie algebra $su(3)$ of  $SU(3)$ is defined as the collection
of $3\times 3$ anti-Hermitian square matrices having trace zero.
The following Gell-Mann matrices, are  the generators of $su(3)$
\[
Z_1=
\begin{bmatrix}
    0 & 1  & 0\\
    1 & 0  & 0\\
    0 & 0  & 0 
\end{bmatrix},
Z_2=i
\begin{bmatrix}
    0 & -1  & 0\\
    1 & 0  & 0\\
    0 & 0  & 0 
\end{bmatrix},
Z_3=
\begin{bmatrix}
    1 & 0  & 0\\
    0 & -1  & 0\\
    0 & 0  & 0 
\end{bmatrix},
Z_4=
\begin{bmatrix}
    0 & 0 & 1\\
    0 & 0 & 0 \\
    1 & 0 & 0 
\end{bmatrix},
\]
\[
Z_5=i
\begin{bmatrix}
    0 & 0  & -1\\
    0 & 0  & 0\\
    1 & 0  & 0 
\end{bmatrix},
Z_6=
\begin{bmatrix}
     0 & 0  & 0\\
    0 & 0  & 1\\
    0 & 1  & 0
\end{bmatrix},
Z_7=i
\begin{bmatrix}
     0 & 0  & 0\\
    0 & 0  & -1\\
    0 & 1  & 0
\end{bmatrix},
Z_8=\frac{1}{\sqrt{3}}
\begin{bmatrix}
     1 & 0  & 0\\
    0 & 1  & 0\\
    0 & 0  & -2
\end{bmatrix}.
\]

The generators $\left\lbrace  Z_2,Z_5,~ \mbox{and}~ Z_7\right\rbrace $ of $su(3)$ are essentially identical to the generators of $so(3)$.  
So, the generators of  $su(3)/so(3)$ are $\left\lbrace  Z_1,Z_3, Z_4, Z_6~ \mbox{and}~ Z_8\right\rbrace $.
So a driftless control system on the symmetric space $SU(3)/SO(3)$ is given by 
\begin{equation}
 \label{random2}
  \dot{x}= x\left(  Z_1u_1+Z_3u_3+Z_4u_4+ Z_6u_6+ Z_8u_8\right) , ~x \in SU(3)/SO(3), ~u_i, i=1,3,4,6,8 \in \mathbb{R}.
  \end{equation}
Since, $[[Z_1,Z_3],Z_4]=  2Z_6$, so, the four generators $\left\lbrace Z_1,Z_3,Z_4, Z_8 \right\rbrace $ upon double bracket commutation generates the full set of basis elements of $su(3)/so(3)$. So the system (\ref{random2}) can be controlled with fewer number of controls and the system can be given by
\begin{equation}
 \label{random3}
  \dot{x}= x\left(  Z_1u_1+Z_3u_3+Z_4u_4+Z_8u_8\right)  , ~x \in SU(3)/SO(3), ~u_i, i=1,3,4,8 \in \mathbb{R}.
  \end{equation}
  
  \markboth{ARCHANA TIWARI}{\rm CONTROL SYSTEM ON A SYMMETRIC SPACE}

\section{Concluding remarks}
In this article, we have extended Rashvesky-Chow's theorem  for global controllability condition on Lie groups to that of symmetric spaces and shown that how it  can be readily and directly implemented for driftless nonlinear control systems on various physical system of symmetric submanifold of $SE(3)$ and random matrix ensembles. Now it is quite interesting to extend such type of  studies for Hamiltonian formalism \cite{Blo1} (using Lie double brackets) and using various analytical and numerical integrators \cite{Cra, Rem} to understand the system dynamics. These type of studies can also be used to understand the controllability of different other physical systems such as quantum systems, molecular systems, etc. \cite{Alt, Ram, Ram1}.

\section*{Acknowledgements} One of the author Prof. K.C. Pati likes to thank SERB, DST for financial support through grant no. EMR/2016/006969.

\end{document}